\documentclass[10pt,draft]{article}
\usepackage{amssymb,amsmath,stmaryrd}
\usepackage[all]{xy}
\newtheorem{theo}{Th\'eor\`eme}[section]
\newtheorem{prop}[theo]{Proposition}
\newtheorem{cor}[theo]{Corollaire}
\newtheorem{lem}[theo]{Lemme}
\newtheorem{defi}[theo]{D\'efinition}

\def \demdu#1 { {\sl D\'emonstration #1.} }
\def \demde#1 { {\sl D\'emonstration de #1 .} } 

\def \Gal {\text{\rm Gal}}

\def \dem{{\sl D\'emonstration.} }

\def \limpro{\lim\limits_{\leftarrow} }

\def \qed{\text {$\quad \square$}}

\def \C {\overline {C}}
\def \U {\overline {U}}

\def \lbigoplus#1#2 { \overset {#2} {\underset {#1} \oplus } }

\def\CM{{\mathbb{C}}}

\def\NM{{\mathbb{N}}}

\def\QM{{\mathbb{Q}}}

\def\ZM{{\mathbb{Z}}}

\def\QG{{\mathfrak Q}}

\def\SG{{\mathfrak S}}

\def\ga{\gamma}
\def\Ga{\Gamma}
\def\Ga{\Gamma}

\def\ep{\varepsilon}

\def\La{\Lambda}

\def\si{\sigma}

\def\ze{\zeta}

\def\GC{{\mathcal{G}}}

\def\Hha{{\widehat{H}}}

\def\Gal{\mathop{\mathrm{Gal}}\nolimits}

\def\Im{\mathop{\mathrm{Im}}\nolimits}

\def\ker{\mathop{\mathrm{ker}}\nolimits}

\def\Tr{\mathop{\mathrm{Tr}}\nolimits}

\def\pr{\prime}

\def \tr {\text {\rm Tr}}
\def\Dbar{\leavevmode\lower.6ex\hbox to 0pt{\hskip-.23ex \accent"16\hss}D}
  \def\cfac#1{\ifmmode\setbox7\hbox{$\accent"5E#1$}\else
  \setbox7\hbox{\accent"5E#1}\penalty 10000\relax\fi\raise 1\ht7
  \hbox{\lower1.15ex\hbox to 1\wd7{\hss\accent"13\hss}}\penalty 10000
  \hskip-1\wd7\penalty 10000\box7}
  \def\cftil#1{\ifmmode\setbox7\hbox{$\accent"5E#1$}\else
  \setbox7\hbox{\accent"5E#1}\penalty 10000\relax\fi\raise 1\ht7
  \hbox{\lower1.15ex\hbox to 1\wd7{\hss\accent"7E\hss}}\penalty 10000
  \hskip-1\wd7\penalty 10000\box7}

\def \cond {\mathrm {cond}}
\title{Sous-modules d'unit\'es en th\'eorie d'Iwasawa}
\author{Jean-Robert Belliard\\
\ \\
\small{\it Universit\'e de Franche-Comt\'e,
Laboratoire de math\'ematiques UMR 6623,} \\
\small{\it 16 route de Gray,
25030 Besan\c con cedex,
France.}\\ \ \\
\small{\texttt{belliard@math.univ-fcomte.fr}}}
\begin{document}
\date{10 septembre 2002}
\maketitle
\setcounter{section}{-1}
\begin{abstract} On donne une condition n\'ecessaire et suffisante,
en termes de ``descente galoisienne'', pour que
le module d'Iwasawa obtenu \`a partir des unit\'es
circulaires \`a la Sinnott soit libre. On d\'etaille ensuite des
exemples qui ne satisfont pas \`a cette condition.
\end{abstract}
\section{Introduction}
Soit $p$ un nombre premier impair fix\'e.
Soit $K$ un corps de nombres et $K_\infty/K$ une $\ZM_p$-extension.
On utilise les notations habituelles qui suivent :

\noindent $\Ga:=\Gal(K_\infty/K)$, $\Ga_n:=\Ga^{p^n}$, $G_n:=\Ga/\Ga_n$ et
 $K_n:=K_\infty^{\Ga_n}$.
Dans toute la suite $\U_n$ d\'esignera le pro-$p$-compl\'et\'e du
quotient des unit\'es de $K_n$ par les racines de l'unit\'e de
$K_n$ et $\La=\ZM_p[[\Ga]]$ l'alg\`ebre d'Iwasawa usuelle. L'objet
de cet article est de d\'egager une condition n\'ecessaire et
suffisante \`a la libert\'e du $\La$-module $\C_\infty$ obtenu par
passage \`a la limite projective \`a partir du
pro-$p$-compl\'et\'e modulo sa torsion du module des unit\'es
circulaires de Sinnott (\cite{Si80}, la d\'efinition est rappel\'ee
en \ref{circsin} ci-dessous). Ce
module est d\'efini lorsque $K$ est ab\'elien sur $\QM$ et c'est  un
sous-module de $\U_\infty:=\limpro \U_n$.\par Dans la premi\`ere
section on redonne une d\'emonstration de crit\`eres g\'en\'eraux
de $\La$-libert\'e, vraisemblablement connus. La proposition
\ref{galdes} fournit un crit\`ere suffisamment fin pour conduire
\`a une \'equivalence si on l'applique au module $\C_\infty$.
Cette \'equivalence (th\'eor\`eme \ref{equiv}) est le r\'esultat
principal de ce travail, elle   est d\'emontr\'ee  dans la
deuxi\`eme section.
 La troisi\`eme partie illustre le th\'eor\`eme
\ref{equiv}, par des contre-exemples \`a la $\La$-libert\'e de
$\C_\infty$. Les cas de $\La$-libert\'e de ce module sont nombreux
et connus (on peut consulter par exemple le \S 1 de \cite{Ku96},
ou \cite {JNT1}). On donne ici une famille infinie d'exemples pour
lesquels on d\'emontre que cette libert\'e n'a pas lieu. D'autres
exemples avec la m\^eme propri\'et\'e ont \'et\'e \'etudi\'es
ind\'ependamment par R. Ku\v cera (\cite{Kuc02}) dans le but
d'\'etablir la diff\'erence entre le module des unit\'es
cyclotomiques de Sinnott et celui de Washington. En cours de route
on \'etablit ici que ces modules sont \'egaux si et seulement si
$\C_\infty$ est $\La$-libre (proposition \ref{compwas}).

\section{G\'en\'eralit\'es}

\subsection{Crit\`eres de $\La$-libert\'e}\label{gencri}

On commence
par rappeler deux lemmes tout \`a fait classiques, dont on redonne une
d\'emonstration pour la commodit\'e du lecteur~:
\begin{lem}\label{coinv} Soit $X$ un $\La$-module de type fini.
Notons $X^\Ga$ (resp. $X_\Ga$)
les invariants (resp. co-invariant) de $X$ sous
l'action de $\Ga$. Alors $X$ est $\La$-libre si et seulement
si $X^\Ga=0$ et $X_\Ga$ est $\ZM_p$-libre.
Dans ce cas la dimension sur $\La$ de $X$ est \'egale \`a celle de $X_\Ga$
sur $\ZM_p$.
\end{lem}
\dem L'implication directe est imm\'ediate. on v\'erifie la r\'eciproque.
Consid\'erons une $\ZM_p$-base $\{\overline {x}_1,\cdots,\overline
{x}_n\}$ de $X_\Ga$ et relevons la en $\{x_i,1\leq i\leq
n\}\subset X$. Par le lemme de Nakayama on obtient d\'ej\`a des
g\'en\'erateurs en nombre voulu. De plus toute relation
non-triviale entre les $x_i$ conduit en premier lieu \`a une
relation non divisible par $T$ entre ces \'el\'ements puisque $X$
est suppos\'e sans $T$-torsion, puis \`a une relation non-triviale
entre les $\overline {x}_i$.\par\qed\par
\begin{lem}\label{indfini}
Soit $Y$ un $\La$-module libre de type fini et $X$ un sous-module
d'indice fini de $Y$. Alors $X$ est $\La$-libre si et seulement si
$X=Y$.\end{lem}
\dem Puisque $X\subset Y$ on a d\'ej\`a $X^\Ga=0$.
Ainsi par le lemme \ref{coinv}
$X$ est $\La$-libre si et seulement si $X_\Ga$ est $\ZM_p$ libre.
La suite de $\Ga$-cohomologie provenant de la suite courte
\xymatrix{0\ar[r]& X\ar[r] & Y \ar[r] & Y/X \ar[r] & 0} fournit
la suite exacte \xymatrix{0\ar[r] & (Y/X)^\Ga \ar[r] & X_\Ga \ar[r] &
Y_\Ga }. Donc $X_\Ga$ est $\ZM_p$-libre si et seulement si $(Y/X)^\Ga$
est sans $\ZM_p$-torsion c'est-\`a-dire trivial.
Comme $Y/X$ est suppos\'e fini on a les
\'equivalences $Y/X=0\iff (Y/X)_\Ga=0 \iff (Y/X)^\Ga =0$.\par\qed\par
Avant de donner des applications, on va \'enoncer un autre
crit\`ere, un peu moins connu. Pour ce, d\'etaillons le contexte
g\'en\'eral dans lequel ce crit\`ere s'applique, contexte qui se
produit tr\`es souvent en arithm\'etique. Il s'agit de la donn\'ee
de deux suites, disons $M_n\subset L_n$, de $\ZM_p[G_n]$-modules
munis, pour $m\geq n$, d'homomorphismes \'equivariants de {\it norme} $N_{m,n}
\colon L_m\longrightarrow L_n$ et d'{\it extension} $i_{n,m}
\colon L_n\longrightarrow L_m$ tels que~:
\begin{itemize}
\item[{1)}] Les restrictions des $N_{n,m}$ et des $i_{m,n}$
d\'efinissent des homomorphismes $N_{m,n}
\colon M_m\longrightarrow M_n$ et $i_{n,m}
\colon M_n\longrightarrow M_m$
\item[{2)}] les compos\'ees $i_{n,m}\circ
N_{m,n}$  co\"\i ncident avec la multiplication par
la {\it trace alg\'ebrique} $\tr_{m,n}:=\sum_{g\in\Gal(K_m/K_n)}
g\in \ZM_p[G_m]$.
\item[{3)}] les compos\'ees $N_{m,n}\circ i_{n,m}$  co\"\i ncident avec  la
multiplication par $p^{m-n}$.
\end{itemize}
Dans ce contexte on a le crit\`ere :
\begin{prop}\label{galdes}
On suppose que les suites $L_n$ et $M_n$ v\'erifient les
conditions suivantes~:
\begin{itemize}
\item[$(i)$] $L_\infty:=\limpro L_n$ est $\La$-libre.
\item[$(ii)$] Les homomorphismes $i_{n,m} \colon L_n\longrightarrow
L_m^{\Gal(K_m/K_n)}$ sont injectifs.
\item[$(iii)$] La suite $M_n$ v\'erifie asymptotiquement la
{\rm "descente galoisienne"}, (i.e. il existe un $N\in \NM$ tel que
pour tout $n\geq N$ on
l'\'egalit\'e $M_{n+1}^{\Gal(K_{n+1}/K_n)}=i_{n,n+1}(M_n)$).
\end{itemize}
Alors $M_\infty:=\limpro M_n$ est aussi $\La$-libre.
\end{prop}
\dem On va appliquer le lemme \ref{coinv} \`a $X=M_\infty$. Les
inclusions $M_n\subset L_n$ donnent par passage \`a la limite
l'inclusion $M_\infty^\Ga\subset L_\infty^\Ga$. Puisque $L_\infty$
est suppos\'e $\La$-libre on obtient donc $M_\infty^\Ga=0$. Ainsi
il reste \`a montrer la libert\'e du $\ZM_p$-module
$(M_\infty)_\Ga$. En appliquant la suite longue de
$\Ga$-cohomologie \`a la suite \xymatrix{ 0\ar[r] & M_\infty
\ar[r]& L_\infty \ar[r] &L_\infty/M_\infty \ar[r] & 0} on
obtient~:

\centerline {\xymatrix {0 \ar[r] & (L_\infty/M_\infty)^\Ga \ar[r]
& (M_\infty)_\Ga \ar[r] & (L_\infty)_\Ga}.}

\noindent La $\La$-libert\'e de $L_\infty$ fournit la
$\ZM_p$-libert\'e de $(L_\infty)_\Ga$, et il suffit donc de
montrer que $(L_\infty/M_\infty)^\Ga$ est sans $\ZM_p$-torsion ou
encore que le sous-$\La$-module fini maximal de
$L_\infty/M_\infty$ est trivial. Notons $F_\infty\subset
L_\infty/M_\infty $ ce sous-module et posons $F_n\subset L_n/M_n$
pour ses projections naturelles. Par compacit\'e de $F_\infty$ on a
$F_\infty=\limpro F_n$, et pour conclure la preuve de \ref{galdes}
il nous reste \`a v\'erifier le~:
\begin{lem}\label{maxf0}
$\forall n\in \NM,\ F_n=0$.\end{lem} \demde{\ref{maxf0}} \ Soit un
$N$ v\'erifiant \ref{galdes} $(iii)$. On
remarque d'abord que la suite des $L_n/M_n$, et donc la suite des
$F_n$, est encore munie d'homomorphismes de norme et d'extension que
l'on note encore $N_{m,n}$ et $i_{n,m}$ par abus de langage. Puisque
$F_\infty$ est fini il existe un $N^\pr\in \NM$
 tel que $\Ga_{N^\pr}$ agit
trivialement sur $F_\infty$. On peut prendre $N^\pr\geq N$.
Puisque par d\'efinition $F_\infty
\longrightarrow F_n$ est surjective $\Ga_{N^\pr}$ agit aussi
trivialement sur tous les $F_n$. De plus il suffit de v\'erifier le lemme
pour $n\geq N^\pr$. Les ordres des $F_n$ sont born\'es, disons par $\#
F_\infty=p^b$.
Fixons donc $n\geq N^\pr$ et soit $f_n\in F_n$. Il existe
alors un $f_{n+b}\in F_{n+b}$ tel que $f_n=N_{n+b,n}(f_{n+b})$ et
donc $i_{n,n+b}(f_n)=i_{n,n+b}\circ N_{n+b,n}(f_{n+b})=\tr_{n+b,n}
(f_{n+b})=p^b f_{n+b}=0$. Ainsi $F_n$ est contenu
dans le noyau : $F_n\subset\ker\left\{ i_{n,n+b} \colon L_n/M_n
\longrightarrow L_{n+b} /M_{n+b}\right\}$.
En appliquant \ref{galdes} $(ii)$ et le lemme du serpent dans le
diagramme

\centerline{\xymatrix{ 0 \ar[r] & M_n \ar[r]
\ar[d]^-{i_{n,n+b}} & L_n  \ar[d]^-{i_{n,n+b}}\ar[r] & L_n/M_n
 \ar[d]^-{i_{n,n+b}} \ar[r] & 0 \\ 0 \ar[r] & M_{n+b} \ar[r] & L_{n+b}
\ar[r] & L_{n+b}/M_{n+b} \ar[r] & 0}}
\noindent on obtient la suite exacte

\centerline{\xymatrix{ 0 \ar[r] & \ker
\left\{ L_n/M_n \rightarrow L_{n+b} /M_{n+b}\right\} \ar[r] & \ker
\left\{ M_{n+b}/M_n \rightarrow L_{n+b}/L_n \right\} }.}
\noindent Ce second noyau est trivial puisque par \ref{galdes} $(iii)$ on a
$$\ker \left\{ M_{n+b}/M_n \rightarrow L_{n+b}/L_n \right\} =\frac
{L_n\cap M_{n+b}}{M_n}\subset
\frac {M_{n+b}^{\Gal(K_{n+b}/K_n)}} { M_n} =0.$$
Cela montre le lemme \ref{maxf0} et la proposition \ref{galdes} suit.\par\qed

\subsection{Applications}

On va appliquer la proposition \ref{galdes} \`a divers
sous-modules d'unit\'es. Lorsqu'elles ont lieu, les propri\'et\'es
$(ii)$ et $(iii)$ de \ref{galdes} sont plut\^ot ais\'ees \`a
v\'erifier. On a donc besoin de trouver un module $\La$-libre qui
puisse jouer le r\^ole de $L_\infty$. Notons $\U^\pr_n$ le
pro-$p$-compl\'et\'e du module des $(p)$-unit\'es de $K_n$ (i.e.
les nombres de valuations triviales en toute place ne divisant pas
$p$) modulo sa torsion (i.e. modulo les racines de l'unit\'es de
$K_n$). Sous certaines hypoth\`eses, un r\'esultat de Greither
montre que $L_n=\U^\pr_n$ convient. R\'esumons ce r\'esultat :
\begin{theo}[\cite{gre94}]
\label{punlib} On suppose ou bien que $K_\infty/K$ est
cyclotomique ou bien que $K$ ne contient pas de racine $p$-i\`eme
de l'unit\'e. Alors $\U^\pr_\infty:=\limpro \U^\pr_n$ est
$\La$-libre de rang $r_1+r_2+g$, o\`u $(r_1,r_2)$ est la signature
de $K$ et $g$ est le nombre de places divisant $p$ qui sont
totalement d\'ecompos\'ees dans $K_\infty/K$ .
\end{theo}
\dem Ce sont le th\'eor\`eme et la proposition 1 de \cite{gre94}, et
lorsque $K_\infty/K$ est cyclotomique c'est un th\'eor\`eme de
Kuz$'$Min (voir \cite{Ku72} theorem 7.2) \par\qed\par
\begin{cor}\label{ulib} Sous les m\^emes hypoth\`eses que \ref{punlib},
 $\U_\infty:=\limpro \U_n$ est $\La$-libre.
\end{cor}
 \dem Il suffit d'appliquer la proposition \ref{galdes} aux
 modules $M_n=U_n$ et $L_n=U^\pr_n$.\par\qed\par

\noindent Signalons que \ref{galdes} en conjonction avec la
conjecture de Leopoldt et un th\'eor\`eme d'Iwasawa (voir
\cite{I73}, theorem 25) permettent de retrouver un cas particulier
du th\'eor\`eme \ref{punlib} en prenant pour modules $M_n=\U^\pr_n$
et pour modu\-les $L_n$ la somme directe des pro-$p$-compl\'et\'es
des groupes multiplicatifs locaux aux places divisant $p$.
\section{Unit\'es cyclotomiques}
Philosophiquement les r\'esultats qui suivent s'appliquent \`a la
distribution d'Iwa\-sawa (voir \cite{BO}, \S 3). Ils s'appliquent
donc aux unit\'es de Stark des extensions ab\'eliennes relatives,
d\`es que ces derni\`eres existent. Ici on se restreint aux
extensions ab\'eliennes sur $\QM$ afin de disposer des (nombreuses
versions) des unit\'es cyclotomiques. On suppose dor\'enavant que
$K$ est un corps de nombres {\it absolument ab\'elien}, et on note
$\GC$ son groupe de Galois. Puisque $p$ est impair on peut aussi
supposer, sans perte de g\'en\'eralit\'e, que $K$ est {\it
totalement r\'eel}. Dans ce contexte, le corollaire
 \ref{ulib} s'applique et fournit la $\La$-libert\'e de
$\U_\infty$. Cela nous permet d'obtenir, pour le module
des unit\'es circulaires de Sinnott une r\'eciproque
au crit\`ere de la proposition \ref{galdes}.
Rappelons la d\'efinition de ce module :
\begin{defi}[\cite{Si80}]
\label{circsin}
Pour $m\in\NM$ on note $\zeta_m=\exp(2i\pi/m)$ la racine
primitive $m^{\text{i\`eme}}$ de l'unit\'e d\'efinie
par le choix d'un plongement de $\overline \QM$ dans $\CM$.
Soit $D$ le sous-$\ZM[\GC]$-module de $K^\times$
multiplicativement engendr\'e par $-1$ et les nombres de la forme
$\ep_{K,m,a}=N_{\QM(\zeta_m),K\cap\QM(\zeta_m)} (1-\ze_m^a)$
o\`u $m$ parcourt $\NM$ et $m\nmid a$. Sinnott d\'efinit
le module des {\rm unit\'es circulaires} de $K$, not\'e $C_K$, comme
l'intersection de $D$ avec les unit\'es $U_K$ de $K$ : $C_K=D\cap
U_K$. Dans la suite on utilisera les notations suivantes :
\begin{center}
 $\C_K:=C_K\otimes \ZM_p$, \hspace*{.7cm}
 $\C_n:=\C_{K_n}$ pour $n\in \NM$,  \hspace*{.7cm}
$\C_\infty:=\limpro \C_n$, \\
 $\widetilde{C}_n:=\Im({\C_\infty \overset {\text{nat}}
\longrightarrow \C_n})$,\hspace*{.7cm}
$\ep_{K,m}=N_{\QM(\zeta_m),K\cap\QM(\zeta_m)}
(1-\ze_m)$ pour $m>1$\end{center}
\end{defi}
Un argument de compacit\'e standard donne l'\'egalit\'e
$\widetilde{C}_n=\bigcap_{m\geq n} N_{K_m,K_n}(\C_m)$. Les suites
de modules $M_n=\C_n$ et $L_n=\U^\pr_n$ sont munies des
homomorphismes de normes et d'extensions habituels qui v\'erifient
trivialement les conditions $1)$, $2)$ et $3)$ du \S \ref{gencri}
et la condition $(ii)$ de la proposition \ref{galdes}. D'apr\`es
le th\'eor\`eme \ref{punlib}, $L_\infty=\U^\pr_\infty$ v\'erifie
l'hypoth\`ese $(i)$ de la m\^eme proposition. On va d\'emontrer
l'\'equivalence suivante~:
\begin{theo}\label{equiv}
Le $\La$-module $\C_\infty$ est libre si et seulement si
la suite $M_n=\C_n$ v\'erifie la condition $(iii)$ de
\ref{galdes}, \`a savoir~:
$$\C_\infty\ \text{est}\ \La\ \text{libre}\iff \exists
\, N\in\NM, \, \forall n\geq N,\,
\C_{n+1}^{\Gal(K_{n+1}/K_n)}=i_{n,n+1}(\C_n)\cong\C_n.$$
\end{theo}
\begin{itemize}
\item[]{\it Remarque:} Dans tous les cas on sait que le $\La$-rang
du module $\C_\infty$ est \'egal \`a celui de $\U_\infty$,
c'est-\`a-dire \'egal \`a $r_1=[K:\QM]$, puisque $K$ est totalement r\'eel.
\end{itemize}
\dem La proposition \ref{galdes} donne l'implication indirecte, on
v\'erifie le sens direct. On suppose donc que $\C_\infty$ est
libre et on va montrer que la suite $\C_n$ v\'erifie la
propri\'et\'e $(iii)$ de la proposition \ref{galdes} qu'on
appellera dor\'enavant la ``{\it descente asymptotique}''. La
premi\`ere \'etape du raisonnement consiste \`a le montrer pour la
suite de modules $\widetilde{C}_n$. Pour ce on utilise un lemme
d\^u \`a Kuz$'$Min~:
\begin{lem}[\cite{Ku72}, 7.3]\label{puinj}
Pour tout $n\in\NM$, l'application naturelle

\centerline {$(\U^\prime_\infty)_{\Ga_n}
\longrightarrow \U^\prime_n$}

\noindent est injective.
\end{lem}\par
\qed\par
\begin{lem}\label{nugaldes}
Si $\C_\infty$ est $\La$-libre alors $\widetilde{C}_n$
v\'erifie la descente asymptotique.
\end{lem}
On pose $Q_\infty:=\U^\prime_\infty/\C_\infty$.
Le lemme du serpent donne alors la suite pour tout $n$~:
\xymatrix{0\ar[r] & Q_\infty^{\Ga_n} \ar[r] &
(\C_\infty)_{\Ga_n}\ar[r] \ar[r]& (\U^\pr_\infty)_{\Ga_n} }.
Le lemme \ref{puinj} et le diagramme commutatif qui suit
permettent donc d'identifier $(Q_\infty)^{\Ga_n}$ avec
le noyau de descente : $(Q_\infty)^{\Ga_n} \simeq
\ker\left((\C_\infty)_{\Ga_n}\twoheadrightarrow \widetilde{C}_n\right)$.

\centerline{\xymatrix{0\ar[r]& (Q_\infty)^{\Ga_n} \ar[r] &
(\C_\infty)_{\Ga_n} \ar[r]\ar@{->>}[d]& (\U^\pr_\infty)_{\Ga_n}
\ar@{^{(}->}[d] \\ \ & \ & \widetilde {C}_n \ar@{^{(}->}[r]&  \U^\pr_n\\}}
\noindent Par noeth\'erianit\'e les noyaux $(Q_\infty)^{\Ga_n}$ se
stabilisent et valent pour $n$ grand
$(Q_\infty)^{\Ga_n}=(Q_\infty)^{\Ga_N}$,
o\`u $N$ est un entier fix\'e. Ainsi pour tout $n\geq N$ on obtient
la suite exacte :

\centerline{\xymatrix{ 0\ar[r] & (Q_\infty)^{\Ga_N}\ar[r] &
 (\C_\infty)_{\Ga_n} \ar[r] & \widetilde{C}_n \ar[r] & 0}}

\noindent Par hypoth\`ese $\C_\infty$ est $\La$-libre donc $
(\C_\infty)_{\Ga_n}$ est $\ZM_p[G_n]$-libre. En outre le groupe
$\Gal(K_{n+1}/K_n)$ agit trivialement sur le $\ZM_p$-module libre
$(Q_\infty)^{\Ga_N}$ puisque $n\geq N$. On en d\'eduit la trivialit\'e
du $0^{\text{i\`eme}}$ groupe de cohomologie modifi\'e \`a la Tate

\centerline{$\Hha^0(K_{n+1}/K_n,\widetilde{C}_n)\simeq
H^1(K_{n+1}/K_n,(Q_\infty)^{\Ga_N})=0$.}
\noindent Mais par d\'efinition des module $\widetilde{C}_n$ la norme
$N_{n+1,n} \colon \widetilde{C}_{n+1} \longrightarrow \widetilde{C}_n$
est surjective.
On obtient donc les \'egalit\'es
$$(\widetilde{C}_{n+1})^{\Gal(K_{n+1}/K_n)}=
\tr_{n+1,n} (\widetilde{C}_{n+1})=i_{n+1,n}\circ N_{n+1,n}
(\widetilde{C}_{n+1})=i_{n+1,n}(\widetilde{C}_n). $$
\par\qed\par
\noindent Pour conclure la preuve de \ref{equiv} on doit maintenant
passer de $\widetilde{C}_n$ \`a $\C_n$. Le lemme suivant montre que
la d\'eviation entre ces deux modules est asymptotiquement constante :
\begin{lem}\label{devuncyc} Soit $n\in \NM$ et
soit $I_n$ le corps d'inertie en $p$ de $K_n$
(i.e. le plus grand sous-corps de $K_n$ dans lequel $p$ ne se
ramifie pas). Alors on a
$$\C_n=\C_{I_n}\ \widetilde{C}_n.$$
La suite de corps $I_n$ est manifestement stationnaire et si $I$
d\'esigne le corps d'inertie de $p$ dans $K_\infty$ on a
l'\'egalit\'e pour tout $n$ tel que $I\subset K_n$ :
$$\C_n=\C_I \ \widetilde{C}_n$$
\end{lem}
\dem La seconde \'egalit\'e est une cons\'equence directe de la
premi\`e\-re. L'inclusion $\C_{I_n}\ \widetilde{C}_n\subset \C_n$
est imm\'ediate. On montre l'inclusion r\'eciproque.
Soit $u\in \C_n$. Par d\'efinition $u$ est une unit\'e de la forme :
$$u=\prod_{m>2} \ep_{K_n,m}^{x_m},$$
o\`u  les $x_m\in \ZM_p[\Gal(K_n/\QM)]$ sont presque tous nuls.
On commence par s\'eparer le produit ci-dessus selon la
divisibilit\'e de $m$ par $p$. On obtient :
$$u=\left(\prod_{m>2,p\nmid m} \ep_{K_n,m}^{x_m}\right)
\left(\prod_{p\nmid m^\pr,a>0} \ep_{K_n,p^a m^\pr}^{x_{p^a m^\pr}}\right).$$
Montrons que le terme $\prod \ep_{K_n,m}^{x_m}$
appartient \`a $\C_{I_n}$.
Soit $R$ l'ensemble des places ramifi\'ees dans $I_n/\QM$ (remarquer
que $p\notin R$).
Si $p\nmid m$ alors $p$ n'est pas ramifi\'e dans $\QM(\ze_m)$ donc
les nombres de la forme $\ep_{K_n,m}$ sont des $R$-unit\'es de $I_n$,
tandis que les nombres de la forme $\ep_{K_n,p^a m^\pr}$ sont des
$p$-unit\'es de $K_n$. Comme il ne peut pas y avoir de compensation
entre les valuations au-dessus de $R$ et celle au-dessus de $p$,
pour que $u$ soit une unit\'e on a n\'ecessairement
$\prod \ep_{K_n,m}^{x_m}\in \U_n$ et
$\prod \ep_{K_n,p^a m^\pr}^{x_{p^a m^\pr}}\in \U_n$.
Puisque $\QM(\ze_m)\cap K_n$ est contenu dans $I_n$ d\`es que $p\nmid
m$, on en d\'eduit $\prod_{m>2,p\nmid m} \ep_{K_n,m}^{x_m}\in \C_{I_n}$.\par
Il reste \`a montrer que le terme
$\prod_{m^\pr,a>0} \ep_{K_n,p^a m^\pr}^{x_{p^a m^\pr}}\in \U_n$
appartient en fait \`a $\widetilde{C}_n$, et par r\'ecurrence sur $n$
il suffit de v\'erifier que cette unit\'e est la norme d'une unit\'e
de $K_{n+1}$ de la m\^eme forme. S\'eparons \`a nouveau ce produit
comme suit :
$$\prod_{p\nmid m^\pr,a>0} \ep_{K_n,p^a m^\pr}^{x_{p^a m^\pr}}=
\left(\prod_{a>0} \ep_{K_n,p^a}^{x_{p^a}}\right) \left( \prod_{p\nmid
m^\pr,m^\pr>1,a>0} \ep_{K_n,p^a m^\pr}^{x_{p^a m^\pr}} \right).
$$
Et posons :
$$y:= \left(\prod_{a>0} \ep_{K_n,p^a}^{x_{p^a}}\right)\quad \text{et}
\quad z:= \left( \prod_{p\nmid m^\pr,m^\pr>1,a>0} \ep_{K_n,p^a
m^\pr}^{x_{p^a m^\pr}} \right)$$ Par le lemme 2.1 $(i)$ de
\cite{So92} les nombres cyclotomiques $(1-\ze_m)$ v\'erifient les
relations de distributions suivantes (pour $r>2$, $r\mid s$ et en
notations additives):
\begin{equation}\label{dist}
\Tr_{\QM(\ze_s),\QM(\ze_r)}(1-\ze_s)=\left(\prod_{l\mid s,l\nmid r}
(1-\si_l^{-1})\right) (1-\ze_r)\end{equation}
o\`u $l$ parcourt l'ensemble des nombres premiers divisant $s$ et pas
$r$ et o\`u $\si_l$ d\'esigne le Frobenius en $l$ de $\Gal(\QM(\ze_s)/\QM)$.
\'Ecrivons $\cond (K_n)=p^c k$ avec $p\nmid k$.
Alors par l'\'equation
(\ref{dist}),
$y$ est une puissance galoisienne de l'uniformisante
cyclotomique $\ep_{K_n,p^c}$ tandis que $z$ est une unit\'e. Donc $y$ est
aussi une unit\'e et peut s'\'ecrire $y=\ep_{K_n,p^c}^t$ avec
$t$ dans l'id\'eal d'augmentation de l'anneau $\ZM_p[\Gal(\QM(\ze_{p^c})\cap
K_n / \QM)]$.
Comme l'\'equation (\ref{dist}) donne aussi l'identit\'e
$\ep_{K_n,p^c}=N_{K_{n+1},K_n}(\ep_{K_{n+1},p^{c+1}})$, on a bien
$y\in \widetilde {C}_n$. \par
 Pour conclure v\'erifions que chaque terme de la forme $\ep_{K_n,p^a
m^\pr}$ avec $a>0$ et $m^\pr>1$ est norme d'unit\'e circulaire de $K_{n+1}$.
Soit $F=\QM(\ze_{p^a m^\pr}) \cap K_n$, $F_1$ le premier \'etage de
la $\ZM_p$-tour sur $F$, et soient $f$  et $f_1$ leurs conducteurs
respectifs. Par l'\'equation
(\ref{dist}) on voit que $\ep_{K_n,p^a
m^\pr}$ est une puissance galoisienne de $\ep_{F,f}$ si $p\mid f$ ou
de $\ep_{F,f}^{1-\sigma_p^{-1}}$ si $p\nmid f$.
Comme en outre seules les places au-dessus de $p$ sont
\'eventuellement ramifi\'ees dans l'extension $F_1/F$ la m\^eme
\'equation donne aussi :
$$N_{F_1/F} (\ep_{F_1,f_1}) = \left\{ \aligned &\ep_{F,f}\quad \text
{si}\quad
p\mid f \\ & \ep_{F,f}^{1-\si_p^{-1}}\quad  \text{si}\quad
p\nmid f \endaligned
\right . $$
Si $F_1$ n'est pas contenu dans $K_n$ alors on a $F_1\subset K_{n+1}$, puis
$N_{K_{n+1},K_n} (\ep_{F_1,f_1})=N_{F_1,F} (\ep_{F_1,f_1})$ et la
preuve est termin\'ee. Sinon comme $\widetilde {C}_n$ est un module
galoisien il suffit de montrer que $\ep_{F_1,f_1}\in \widetilde{C}_n$
ce qui conduit \`a terme au cas pr\'ec\'edent.
%Cela se voit sur un syst\`eme de g\'en\'erateur de
%$\C_n$ et avec les relations de distributions v\'erifi\'ees par les
%nombres cyclotomiques. Soit $R$ l'ensemble des places ramifi\'ees dans
%$I$,
%alors les nombres de la forme $\ep_{K_n,p^af}$ sont
%des normes universelles  de $p$-unit\'es
%(donc des \'el\'ements de $\widetilde{C}_n$
%pour $f>1$) tandis que les nombres $\ep_{K_n,f}$  sont des $R$-unit\'es
%de $I_n$ d\`es que $p\nmid f$.
%Il faut ensuite prendre l'intersection des nombres
%cyclotomiques avec les unit\'es, ce qui ne pose pas de probl\`eme.
%Pour les d\'etails voir
%\cite{JNT1} lemme 1.3 et proposition 1.6.
\par\qed\par
\noindent On reprend la preuve du th\'eor\`eme \ref{equiv}. Fixons
un $n$ suffisamment grand pour avoir $I\subset K_n$. Alors on a
$\C_I^{\Gal(K_{n+1}/K_n)}=\C_I$ dont on d\'eduit :
$$\begin{aligned} \C_{n+1}^{\Gal(K_{n+1}/K_n)}&=
\left (\C_I\ \widetilde{C}_{n+1}\right )^{\Gal(K_{n+1}/K_n)} =
\C_I \left(\widetilde{C}_{n+1}\right)^{\Gal(K_{n+1}/K_n)}\\ &=
\C_I \widetilde{C}_n=\C_n\ .\end{aligned}$$ Ce qui d\'emontre le
th\'eor\`eme \ref{equiv}.\par\qed\par

\section{Exemples et contre-exemples}
A partir du th\'eor\`eme \ref{equiv} qui donne une
\'equivalence, il est naturel de se demander si les deux
alternatives peuvent effectivement se produire. Lorsque le corps
$K$ est le sous-corps r\'eel maximal d'un corps cyclotomique, on
sait par des r\'esultats sur les distributions \`a la Kubert-Lang
que $\C_\infty$ est $\La$-libre (voir par exemple \cite{Ku96}).
Plus g\'en\'eralement, par un passage \`a la limite imm\'ediat, on
peut d\'eduire de \cite{JNT1} une condition suffisante et d'autres
exemples o\`u $\C_\infty$ est $\La$-libre. D'autre part on sait que les
unit\'es circulaires \`a la Sinnott ne v\'erifient pas toujours la
"descente galoisienne" c'est-\`a-dire qu'il existe des corps
ab\'eliens $F\subset L$ tels que $C_F\subsetneq C_L\bigcap F$
(voir \cite{Grei93}). En particulier \ref{equiv} n' entra\^\i ne
pas la libert\'e de $\C_\infty$ : on va pr\'esenter une liste
d'exemples de corps $K$ pour lesquels cette libert\'e n'a pas
lieu. Pour produire
ces contre-exemples on va utiliser des hypoth\`eses de
d\'ecomposition assez contraignantes dans le style de
\cite{Grei93}. Ces hypoth\`eses permettent de mieux contr\^oler la
structure galoisienne des unit\'es circulaires et simplifient
notablement les calculs qui suivent. Ce fait justifie la
terminologie {\it ``g\"unstige $(p+1)$-tuple''} de \cite{Grei93}.
\'Enon\c cons ces conditions. A partir d'ici, et jusqu\`a la fin
de ce paragraphe, on suppose que le corps de nombres ab\'elien $K$ et
le nombre premier $p$ v\'erifient les trois conditions suivantes :
\begin{itemize}
\item[1--] le conducteur de  $K$   est le produit sans facteurs carr\'es
$f_K=\prod_{i=1}^{i=p+1} l_i$ de $p+1$ nombres premiers $l_i$ tels
que $l_i\equiv 1 [p]$ et pour chaque $j\neq i$ il existe un entier
$x_{i,j}$ tel que $l_i\equiv x_{i,j}^p [l_j]$.
\item[2--] $\GC:=\Gal(K/\QM) \simeq (\ZM/p\ZM)^2$
\item[3--] Notons $K^1,K^2,\ldots,$ $K^{p+1}$ les $(p+1)$ sous-corps
non triviaux de $K$. On suppose en outre que quitte \`a
renum\'eroter ces sous-corps, leur conducteur vaut~:
$$\cond(K^j)=\prod_{i=1,i\neq j}^{i=p+1} l_i$$
\end{itemize}
\begin{itemize}
\item [{}]{\it Remarque~:}
\begin{itemize}
\item [$(i)$] En reprenant  le raisonnement de
\cite{Grei93}
qui s'appuie sur le th\'eor\`eme
de densit\'e de \v Cebotarev  on peut d\'emontrer
l'existence (d'une infinit\'e) de $(p+1)$-uples
de premier $(l_i)_{i=1}^{p+1}$ tels que $p$ et tout sous-corps de
$\QM(\zeta_{\prod l_i})$ v\'erifient la condition 1.
Ce sont de tels $(p+1)$-uples qui furent d\'enomm\'es {\it g\"unstige
$(p+1)$-tuple}  dans \cite{Grei93}.
\item[$(ii)$] Ensuite \'etant fix\'e un tel g\"unstige $(p+1)$-uple
$l_1,\ldots, l_{p+1}$ on peut v\'erifier que $\QM(\ze_{\prod l_i})$
contient un sous-corps $K$ v\'erifiant 2- et 3- : c'est un exercice sur
les groupes ab\'eliens laiss\'e au lecteurs.
\item[$(iii)$] Alternativement l'exemple du \S IV.2 de \cite{JNT1}
v\'erifie les conditions ci-dessus pour $p=3$.
\end{itemize}
\end{itemize}
 La condition 3-- entra\^ine que les Frobenius $\left (\frac {l_i}
{K^i/\QM}\right )$
sont des \'el\'ements bien d\'efinis de
$\Gal(K^i/\QM)=\Gal(K^i_\infty/\QM_\infty)$.
De plus ces Frobenius sont triviaux par la condition 1--~:
\begin{itemize}
\item[4--] $\forall i, 1\leq i\leq p+1 \quad :
\quad \left (\frac {l_i} {K^i/\QM}\right )
=\left (\frac {l_i} {K^i_\infty/\QM_\infty}\right )=1$.
\end{itemize}
Avec toutes ces hypoth\`eses on va d\'emontrer :
\begin{theo}\label{contrex} $\C_\infty$ n'est pas $\La$-libre.
\end{theo}
\dem
On commence par se donner un syst\`eme de g\'en\'erateurs galoisiens
de  $\C_\infty$. Pour tout corps ab\'elien $F$, de conducteur $f$
disons, on note $\ep_F$ le nombre cyclotomique
$$\ep_F =\ep_{F,f}= N_{\QM(\zeta_f),F} (1-\zeta_f)$$
On d\'esignera par $\ep_\infty^i$ (pour $1\leq i\leq p+1$),
 $\ep_\infty^\QM$ et $\ep_\infty^K$
les \'el\'ements suivants de $\C_\infty$ :
$$\aligned
\ep^i_\infty &:= (N_{K^i_1,K^i} (\ep_{K_1^i}), \ep_{K_1^i}, \ldots ,
\ep_{K_n^i},\ldots)_{n\geq 1} \\
\ep^\QM_\infty &:= (\ga-1) (p,\ep_{\QM_1},\ldots,\ep_{\QM_n},\ldots
)_{n\geq 1}\\
\ep^K_\infty &:=
(N_{K_1,K}(\ep_{K_1}),\ep_{K_1},\ldots,\ep_{K_n},\ldots)_{n\geq
1}\endaligned.$$
A partir du lemme 2.3 de \cite{Grei92} on peut v\'erifier
le lemme suivant (remarquer cependant
que les g\'en\'erateurs utilis\'es ici diff\`erent de ceux de \cite{Grei92}
 par des inversibles de $\La[\GC]$) :
\begin{lem} Le syst\`eme $(\ep^I_\infty)_I$, o\`u $I$ parcourt
l'ensemble $\{K,\QM,1,2,\ldots,p+1\}$, engendre $\C_\infty$ sur
$\La[\GC]$.\par\qed\end{lem}
Par un passage \`a la limite projective le long de la $\ZM_p$-tour
sur l'\'equation (\ref{dist})
on obtient les relations (en notation additive) :
$$\begin{aligned} \left (R_1^{(i)}\right ) &: Tr_{K_\infty,K_\infty^i}
(\ep^K_\infty)= \left (1-\left ( \frac {l_i} {K^i_\infty/\QM}
 \right )^{-1} \right ) \ep_\infty^i\\
\left (R_2^{(i)}\right ) &: Tr_{K^i_\infty,\QM_\infty} (\ep^i_\infty)
= \left (\frac 1 {\ga -1} \prod_{j\neq i}
\left (1-\left ( \frac {l_j} {\QM_\infty/\QM} \right )^{-1} \right )\right )
 \ep_\infty^\QM
\end{aligned} $$
O\`u $Tr_{F,L}(x)=(\sum_{g\in\Gal(F/L)} g) x $ d\'esigne la trace alg\'ebrique
d\'efinie pour toute extension galoisienne finie $F/L$. A priori
tous les termes $\left (1-\left ( \frac {l_i} {K^i_\infty/\QM}
 \right )^{-1} \right )$ appartiennent \`a $\La[\Gal(K^i/\QM)]$.
Mais gr\^ace \`a 4--, ces facteurs sont en fait des \'el\'ements de
$\ZM_p[\Ga]\subset\La$.
En combinant n'importe quel couple de relations avec le m\^eme $i$
on obtient :
$$\left (R_3\right ) : Tr_{K_\infty,\QM_\infty} (\ep^K_\infty)
= \left (\frac 1 {\ga -1} \prod_{i=1}^{i=p+1}
\left (1-\left ( \frac {l_j} {\QM_\infty/\QM} \right )^{-1} \right )\right )
 \ep_\infty^\QM.$$
A l'aide de l'identit\'e formelle $p = -Tr_{K_\infty,\QM_\infty} +
\sum_{i=1}^{i=p+1} Tr_{K_\infty,K^i_\infty}$, des relations
$(R^{(i)}_1)$ et $(R_3)$ on obtient une autre relation utile :
$$\aligned
\left (R_4\right ) :
 p (\ep^K_\infty)
&= \left (- \frac 1 {\ga -1}  \prod_{i=1}^{i=p+1}
\left (1-\left ( \frac {l_j} {\QM_\infty/\QM} \right )^{-1} \right )\right )
 \ep_\infty^\QM
\\ &+ \sum_{i=1}^{i=p+1}
\left (1-\left ( \frac {l_i} {K^i_\infty/\QM}
 \right )^{-1} \right ) \ep_\infty^i \endaligned$$
On consid\`ere maintenant le sous-$\La[\GC]$-module de $\C_\infty$,
disons $\mathfrak {S}_\infty$,
engendr\'e par tous les $\ep^I_\infty$ sauf  $\ep^K_\infty$~:
$$\mathfrak {S}_\infty := \left \langle \left \{ \ep^I_\infty, I = \QM, 1,
\ldots, p+1\right \} \right \rangle $$
On note $\mathfrak {Q}_\infty$ le quotient
$$\mathfrak {Q}_\infty:= \C_\infty/\mathfrak {S}_\infty $$
Bien sur $\mathfrak {Q}_\infty$ est $\La[\GC]$-monog\`ene, engendr\'e par
$\ep^K_\infty + \mathfrak {S}_\infty$, et tu\'e par $p$ \`a cause de $(R_4)$.
Il suit que $\SG_\infty$ a le m\^eme $\La$-rang que $\C_\infty$,
\`a savoir $p^2$. Par un choix appropri\'e de $\La$-base on peut voir que
$\SG_\infty$ est $\La$-libre :
\begin{lem}\label{Sfree} Soit $g_i$ un g\'en\'erateur de
$\Gal(K^i/\QM) \cong \Gal
(K^i_\infty/\QM_\infty)\simeq \ZM/p\ZM$. Alors l'ensemble qui suit
forme une $\La$-base du $\La$-module
libre $\SG_\infty$ :
$$\left \{
\begin{array}{rccl} &  & \ep^\QM_\infty,& \\
 \ep^1_\infty,& g_1 \ep^1_\infty, &\ldots ,&
g_1^{p-2} \ep^1_\infty, \\
 \ep^2_\infty,& g_2 \ep^2_\infty, & \ldots,&
g_2^{p-2} \ep^2_\infty,\\
& & \vdots & \\
& & g_i^j \ep_\infty^i, & \\
 & & \vdots &\\
 \ep^{p+1}_\infty,& g_{p+1} \ep^{p+1}_\infty,
& \ldots, & g_{p+1}^{p-2} \ep^{p+1}_\infty \end{array} \right \}$$
\end{lem}
\dem On a choisi $1+(p+1)(p-1)=p^2$ \'el\'ements de $\SG_\infty$
donc il suffit  de v\'erifier que ces \'el\'ements forment un
syst\`eme g\'en\'erateur. En fait \`a partir du syst\`eme
g\'en\'erateur de la d\'efinition de $\SG_\infty$ on a seulement
retir\'e un conjugu\'e de chaque $\ep^i_\infty$, et la relation
$(R_2^i)$ permet d'\'ecrire ce conjugu\'e comme combinaison
lin\'eaire des autres conjugu\'es et de
$\ep^\QM_\infty$.\par\qed\par
\begin{lem}\label{Q<>0} $\QG_\infty \neq \{0\}$ \end{lem}
\dem Cela \'equivaut \`a $\ep^F_\infty \notin \SG_\infty$, qui revient
aussi \`a $ p \ep^F_\infty \notin p \SG_\infty$. Mais la relation $(R_4)$
donne les coefficients de $p \ep^F_\infty$ dans la base de
$\SG_\infty$ du lemme
\ref{Sfree} :
$$\aligned p (\ep^K_\infty)
&= \left (- \frac 1 {\ga -1}  \prod_{i=1}^{i=p+1}
\left (1-\left ( \frac {l_j} {\QM_\infty/\QM} \right )^{-1} \right )\right )
 \ep_\infty^\QM\\
&+ \sum_{i=1}^{i=p+1}
\left (1-\left ( \frac {l_i} {K^i_\infty/\QM}
 \right )^{-1} \right ) \ep_\infty^i\endaligned$$
Or les termes de la forme $\left (1-\left ( \frac {l_i} {K^i_\infty/\QM}
 \right )^{-1} \right )$ appartiennent \`a $\La$ et sont premiers \`a
$p$. \par\qed\par
On est maintenant arm\'e pour donner explicitement un \'el\'ement
d'ordre $p$ dans
$(\C_\infty)_\Ga$.
\begin{lem}\label{torexp}
$\ep^K_\infty + (\ga -1) \C_\infty$ est d'ordre $p$ dans
$(\C_\infty)_\Ga$.\end{lem}
\dem
En usant \`a nouveau de la relation $(R_4)$ on constate que
$p\ep^K_\infty \in (\ga - 1) \C_\infty$.
Il reste donc \`a v\'erifier que $\ep^K_\infty$ lui-m\^eme n'appartient pas
$(\ga - 1) \C_\infty$. Consid\'erons le diagramme :
\par
\noindent \centerline {\xymatrix{\C_\infty
\ar@{->>}[r] \ar@{->>}[d] & \QG_\infty \ar@{->>}[d] \\
(\C_\infty)_\Ga \ar@{->>}[r] & (\QG_\infty)_\Ga \\}}
Il suffit de montrer que $\ep^K_\infty + (\ga -1) \C_\infty$ ne
s'envoie pas sur $0$ dans $(\QG_\infty)_\Ga$.
Mais par d\'efinition de $\QG_\infty$, les conjugu\'es $g\ep^K_\infty$
pour $g\in\GC$ l'engendrent, de sorte que leurs images engendrent
$(\QG_\infty)_\Ga$.  Et comme $\GC$ agit par automorphisme sur
$(\QG_\infty)_\Ga$, les conjugu\'es $g \ep^F_\infty$ s'envoient
simultan\'ement sur $0$ ou non.
Ainsi $\ep^F_\infty \in (\ga-1)
\C_\infty$ contredirait, par le lemme de Nakayama, le lemme \ref{Q<>0}.
Ceci prouve le lemme \ref{torexp}
et le th\'eor\`eme \ref{contrex}. \par\qed\par
Par une autre construction R. Ku\v cera dans
\cite{Kuc02} donne un corps de nombre $K$ pour lequel les unit\'es
circulaires \`a la Sinnott diff\`erent au niveau infini des {\it unit\'es
cyclotomiques \`a la Washington} (avec $p=3$; voir \cite{KN95} et
\cite{Kuc02} pour la terminologie ``unit\'es de Washington'').
Notons $\C^W_n$ le tensoris\'e avec
$\ZM_p$ des unit\'es cyclotomique de $K_n$ \`a la Washington,
et $\C^W_\infty$ leur limite projective.
Par construction et avec le corollaire 1 de \cite{GK89} les
unit\'es de Washington v\'erifient la ``descente galoisienne''.
Ce module $\C^W_\infty$ est donc $\La$-libre par la
proposition \ref{galdes}.
 On a alors :
\begin{prop}\label{compwas} $\C_\infty$ est $\La$-libre si et seulement si
$\C_\infty=\C^W_\infty$.\end{prop}
\dem D'apr\`es le r\'esultat principal de \cite{KN95} le quotient
$\C^W_\infty/\C_\infty$ est fini. La proposition suit donc du lemme
\ref{indfini}.\par \qed\par
La construction de \cite{Kuc02} fournit donc un autre
exemple o\`u $\C_\infty$ n'est pas $\La$-libre. Dans cet exemple
$7=2p+1$ nombres premiers sont ramifi\'es dans $K$. Dans la
famille d'exemples d\'etaill\'es ici, avec $p=3$, on obtient $4$
nombres premiers ramifi\'es. D'apr\`es \cite{JNT1} on sait que si $2$
nombres premiers (au plus) sont ramifi\'es dans $K$ alors $\C_\infty$ est
$\La$-libre. Actuellement dans le cas particulier o\`u exactement $3$
nombres premiers distincts de $p$ sont ramifi\'es dans le corps de
base $K$ on ne dispose ni de preuve de la $\La$-libert\'e de
$\C_\infty$ ni de contre-exemple.

\def\Dbar{\leavevmode\lower.6ex\hbox to 0pt{\hskip-.23ex \accent"16\hss}D}
  \def\cfac#1{\ifmmode\setbox7\hbox{$\accent"5E#1$}\else
  \setbox7\hbox{\accent"5E#1}\penalty 10000\relax\fi\raise 1\ht7
  \hbox{\lower1.15ex\hbox to 1\wd7{\hss\accent"13\hss}}\penalty 10000
  \hskip-1\wd7\penalty 10000\box7}
  \def\cftil#1{\ifmmode\setbox7\hbox{$\accent"5E#1$}\else
  \setbox7\hbox{\accent"5E#1}\penalty 10000\relax\fi\raise 1\ht7
  \hbox{\lower1.15ex\hbox to 1\wd7{\hss\accent"7E\hss}}\penalty 10000
  \hskip-1\wd7\penalty 10000\box7}
\providecommand{\bysame}{\leavevmode\hbox to3em{\hrulefill}\thinspace}

%\nocite{Wa}
%\bibliography{biblio}

\end{document}